\documentclass[11pt,a4paper]{amsart}
\usepackage{graphicx}
\usepackage{amssymb}
\usepackage{amsmath}
\usepackage{amsthm}
\usepackage{amscd}
\usepackage[all,2cell,dvips]{xy}
\UseAllTwocells \SilentMatrices
\newtheorem{thm}{Theorem}[section]

\newtheorem{cor}[thm]{Corollary}

\newtheorem{exm}[thm]{Example}

\theoremstyle{definition}

\theoremstyle{remark}
\newtheorem{rem}[thm]{\bf Remark}
\numberwithin{equation}{section}

\begin{document}
\title[On  extensions of covariantly finite subcategories]
{On  extensions of covariantly finite subcategories}
\author[  Xiao-Wu Chen
] {Xiao-Wu Chen}
\thanks{This project was supported by China Postdoctoral Science Foundation No. 20070420125, and
was also partially supported by the National Natural Science
Foundation of China (Grant No.s 10501041 and 10601052). The author
also gratefully acknowledges the support of K. C. Wong Education
Foundation, Hong Kong}
\thanks{E-mail:
xwchen$\symbol{64}$mail.ustc.edu.cn}
\maketitle
\date{}%
\dedicatory{}%
\commby{}%
\begin{center}
Department of Mathematics\\
 University of Science and
Technology of China \\Hefei 230026, P. R. China
\end{center}

\begin{abstract}
In \cite{GT}, Gentle and Todorov proved  that in an abelian category
with enough projective objects, the extension subcategory of two
covariantly finite subcategories is still covariantly finite.  We
give an counterexample to show that  Gentle-Todorov's theorem may
fail in arbitrary abelian categories;  we also prove that a
triangulated version of  Gentle-Todorov's theorem holds; we make
applications of Gentle-Todorov's theorem to obtain short proofs to a
classical result by Ringel and a recent result by Krause and
Solberg.
\end{abstract}

\section{Main Theorems}

Let $\mathcal{C}$ be an additive category. By a subcategory
$\mathcal{X}$ of $\mathcal{C}$ we always mean a full additive
subcategory. Let $\mathcal{X}$ be a subcategory of $\mathcal{C}$ and
let $M\in\mathcal{C}$. A morphism $x_M: M\longrightarrow X_M$ is
called a \emph{left $\mathcal{X}$-approximation} of $M$ if $X_M\in
\mathcal{X}$ and every morphism from $M$ to an object in
$\mathcal{X}$ factors through $x_M$. The subcategory $\mathcal{X}$
is said to be \emph{covariantly finite} in $\mathcal{C}$, if every
object in $\mathcal{C}$ has a left $\mathcal{X}$-approximation. The
notions of left $\mathcal{X}$-approximation and covariantly finite
are also known as \emph{$\mathcal{X}$-preenvelop} and\emph{
preenveloping}, respectively. For details, see \cite{AS1, AS2} and
\cite{E}. \vskip 5pt

To state our main result, let $\mathcal{C}$ be an abelian category.
Let $\mathcal{X}$ and $\mathcal{Y}$ be its subcategories. Set
$\mathcal{X}\ast \mathcal{Y}$ to be the subcategory consisting of
objects $Z$ such that there is a short exact sequence
$0\longrightarrow X \longrightarrow Z \longrightarrow Y
\longrightarrow 0$ with $X\in \mathcal{X}$ and $Y\in \mathcal{Y}$,
and it is called the \emph{extension subcategory} of $\mathcal{Y}$
by $\mathcal{X}$. Note that the operation ``$\ast$" on subcategories
is associative. Recall that an abelian category $\mathcal{C}$
\emph{has enough projective objects}, if for each object $M$ there
is an epimorphism $P\longrightarrow M$ with $P$ projective.

\vskip 10pt

The following result is due to Gentle and Todorov \cite{GT}, which
extends the corresponding results in artin algebras and coherent
rings, obtained by Sikko and Smal{\o} (see \cite[Theorem 2.6]{SS}
and \cite{SS2}).

\begin{thm} \label{mainthm}
{\rm (Gentle-Todorov \cite[Theorem 1.1, ii)]{GT})}Let $\mathcal{C}$
be an abelian category with enough projective objects. Assume that
both $\mathcal{X}$ and $\mathcal{Y}$ are covariantly finite
subcategories in $\mathcal{C}$. Then the extension subcategory
$\mathcal{X}\ast \mathcal{Y}$ is covariantly finite.\end{thm}

\noindent {\bf Proof.}\quad  The proof is also due to Gentle and
Todorov, and here we just include it for an inspiration of the proof
of Theorem 1.4. Note that the argument here resembles the one in the
proof of \cite[Lemma 1.3]{KS}.

Assume that $M\in \mathcal{C}$ is an arbitrary object. Take its left
$\mathcal{Y}$-approximation $y_M: M \longrightarrow Y_M$ with
$Y_M\in \mathcal{Y}$. By assumption the category $\mathcal{C}$ has
enough projective objects, we may take an epimorphism $\pi_M:
P\longrightarrow Y_M$ with $P$ projective. Consider the short exact
sequence $0\longrightarrow K \longrightarrow M\oplus P
\stackrel{(y_M,\; \pi_M)}\longrightarrow Y_M \longrightarrow 0$.
Take a left $\mathcal{X}$-approximation $x_K: K \longrightarrow X_K$
of $K$. We form a pushout and then get the following commutative
exact diagram

\[\xymatrix@C=18pt @R=7.5pt{
0 \ar[rr] && K \ar[dd]_{x_K} \ar[rr] &&  M\oplus P
\ar@{.>}[dd]^{(z_M,\;
b_M)} \ar[rr]^-{(y_M, \; \pi_M)} && Y_M \ar@{=}[dd] \ar[rr] && 0\\
& &  & (*) & & & \\
 0 \ar[rr] && X_K \ar@{.>}[rr]^-{i_M} && Z_M \ar[rr] && Y_M
\ar[rr] && 0.
 }\]

Note that $Z_M\in \mathcal{X}\ast \mathcal{Y}$. We claim that the
morphism $z_M: M \longrightarrow Z_M$ is a left $\mathcal{X}\ast
\mathcal{Y}$-approximation of $M$. Then we are done.

\vskip 5pt

To see this, assume that we are given a morphism $f: M
\longrightarrow Z$, and $Z\in \mathcal{X}\ast \mathcal{Y}$, that is,
there is an exact sequence $0 \longrightarrow X
\stackrel{i}\longrightarrow Z \stackrel{\pi}\longrightarrow Y
\longrightarrow 0$ with $X\in \mathcal{X}$ and $Y\in \mathcal{Y}$.
Since $y_M$ is a left $\mathcal{Y}$-approximation, then the
composite morphism $\pi\circ f$ factors through $y_M$, say there is
a morphism $c: Y_M\longrightarrow Y $ such that $\pi\circ f=c\circ
y_M$. Consider the composite morphism
$P\stackrel{\pi_M}\longrightarrow Y_M \stackrel{c}\longrightarrow Y$
and the epimorphism $\pi: Z \longrightarrow Y$. Since $P$ is
projective, we have a morphism $b: P \longrightarrow Z$ such that
$\pi\circ b=c\circ \pi_M$. Hence we have the following  commutative
exact diagram.

\[\xymatrix@C=40pt{
0 \ar[r] & K \ar[r] \ar@{.>}[d]_{a}  & M\oplus P \ar[d]^-{(f, \; b)}
\ar[r]^-{(y_M, \; \pi_M)} & Y_M \ar[d]^-{c} \ar[r] & 0\\
0\ar[r] & X \ar[r]^-{i} & Z \ar[r]^-{\pi} & Y \ar[r] & 0. }\]

Since $x_K$ is a left $\mathcal{X}$-approximation, the morphism $a$
factors through $x_K$, say we have $a=a'\circ x_K$ with $a': X_K
\longrightarrow X$. Then by the universal mapping property of the
pushout square $(*)$, there is a unique morphism $h: Z_M
\longrightarrow Z$ such that
\begin{align*}
h\circ i_M=i\circ a', \quad \mbox{and }\; h\circ (z_M,\; b_M)=(f, \;
b).
\end{align*}
In particular, the morphism $f$ factors through $z_M$, as required.
\hfill $\blacksquare$

\vskip 10pt

\begin{rem}
The same proof yields another version of Gentle-Todorov's theorem
(\cite[Theorem 1.1, i)]{GT}): Let $\mathcal{C}$ be an abelian
category, $\mathcal{X}$ and $\mathcal{Y}$ covariantly finite
subcategories in $\mathcal{C}$. Assume further that $\mathcal{Y}$ is
closed under subobjects. Then the extension subcategory
$\mathcal{X}\ast \mathcal{Y}$ is covariantly finite. (In the proof
of this case, $y_M$ could be assumed to be epic, and then we just
take $P=0$.) \hfill $\blacksquare$
\end{rem}

\vskip 10pt

\begin{exm} {\rm We remark that Gentle-Todorov's theorem  may fail if the abelian category has
not enough projective objects. To give an example, let $k$ be a
field, and let $Q$ be the following quiver
\[
\xymatrix@C=50pt @R=30pt { \cdot_1  \ar@(l,u)^{\bar{\alpha}}
\ar[r]^\beta & \cdot_2 } \]
 where $\bar{\alpha}=\{\alpha_i\}_{i\geq
1}$ is  a  family of arrows indexed by positive numbers. Recall that
a representation of $Q$, denoted by $V=(V_1, V_2; V_{\bar{\alpha}},
V_\beta)$, is given by the following data: two $k$-spaces $V_1$ and
$V_2$ attached to the vertexes, and
$V_{\bar{\alpha}}=(V_{\alpha_i})_{i\geq 1}$, $V_{\alpha_i}:
V_1\longrightarrow V_1$ and $V_\beta: V_1\longrightarrow V_2$
$k$-linear maps attached to the arrows. A morphism of
representations, denoted by $f=(f_1, f_2): V\longrightarrow V'$,
consists of two linear maps $f_i: V_i\longrightarrow V'_i$, $i=1,2$,
which are compatible with the linear maps attached to the arrows.
Denote by $\mathcal{C}$ the category of representations $V=(V_1,
V_2; V_{\bar{\alpha}}, V_\beta)$ of $Q$ such that ${\rm dim}_k\;
V={\rm dim}_k \; V_1+{\rm dim}_k \; V_2<\infty$ and $V_{\alpha_i}$
are zero for all but finitely many $i$'s. Then $\mathcal{C}$ is an
abelian category with finite-dimensional Hom spaces. \vskip 5pt

 Denote by $S_i$ the
one-dimensional representation of $Q$ seated on the vertex $i$ with
zero maps attached to all the arrows, $i=1,2$. Consider the
two-dimensional representation $M=(M_1=k, M_2=k; M_{\bar{\alpha}}=0,
M_\beta=1)$. Denote by $\mathcal{X}$ (resp. $\mathcal{Y}$) the
subcategory consisting of direct sums of copies of $S_1$ (resp.
$M$). Then both $\mathcal{X}$ and $\mathcal{Y}$ are covariantly
finite in $\mathcal{C}$. However we claim that
$\mathcal{Z}=\mathcal{X}\ast \mathcal{Y}$ is not covariantly finite.

\vskip 5pt

 In fact, the
representation $S_2$ does not have a left
$\mathcal{Z}$-approximation. Otherwise, assume that  $\phi:
S_2\longrightarrow V$ is a left $\mathcal{Z}$-approximation. Take
$i_0>>1$ such that $V_{\alpha_{i_0}}$ is a zero map. Consider the
following three-dimensional representation
$$W=(W_1=\begin{pmatrix} k \\ k\end{pmatrix}, W_2=k;
W_{\alpha_{i_0}}=\begin{pmatrix} 0 & 0 \\
                                1 & 0\end{pmatrix}, W_{\alpha_i}=0 \mbox{ for } i\neq i_0, W_\beta=(1, \; 0)).$$
We have a non-split exact sequence of representations
$0\longrightarrow S_1 \longrightarrow W\longrightarrow
M\longrightarrow 0$. In particular, $W\in \mathcal{Z}$. Note that
${\rm Hom}_\mathcal{C}(S_2, W)\simeq k$. Hence there is a morphism
$f=(f_1, f_2): V\longrightarrow W$ such that $f\circ \phi\neq 0$.
However this is not possible.  Note that
$W_{\alpha_{i_0}}(f_1(V_1))=f_1(V_{\alpha_{i_0}}(V_1))=0$ by the
choice of $i_0$, and that ${\rm Ker}\; W_{\alpha_{i_0}}={\rm Ker}\;
W_\beta$, we obtain that $W_\beta(f_1(V_1))=0$. Note that both the
representations $S_1$ and $M$ satisfy that the map attached to the
arrow $\beta$ is surjective, and then by Snake Lemma we infer that
every representation in $\mathcal{Z}$  has this property, in
particular, the representation $V$ has this property, that is,
$V_2=V_\beta(V_1)$. Hence we have
$0=W_\beta(f_1(V_1))=f_2(V_\beta(V_1))=f_2(V_2)$, and thus we deduce
that $f_2=0$. This will force that the composite
$S_2\stackrel{\phi}\longrightarrow V\stackrel{f}\longrightarrow W$
is zero. \hfill $\blacksquare$ }
\end{exm}

\vskip 10pt

We also have  a triangulated version of Gentle-Todorov's theorem.
Let $\mathcal{C}$ be a triangulated category with the translation
functor denoted by $[1]$. For triangulated categories, we refer to
\cite{V1, K}. Let $\mathcal{X}$, $\mathcal{Y}$ be its subcategories.
Set $\mathcal{X}\ast \mathcal{Y}$ to be the \emph{extension
subcategory}, that is, the subcategory consisting of objects $Z$
such that there is a triangle $X \longrightarrow Z \longrightarrow Y
\longrightarrow X[1]$ with $X\in \mathcal{X}$ and $Y\in
\mathcal{Y}$. Again this operation ``$\ast$" on subcategories is
associative by the octahedral axiom (TR4). Then we have the
following result.

\begin{thm} \label{mainthmtri} Let $\mathcal{C}$ be a triangulated category.
Assume that both $\mathcal{X}$ and $\mathcal{Y}$ are covariantly
finite subcategories in $\mathcal{C}$. Then the extension
subcategory $\mathcal{X}\ast \mathcal{Y}$ is covariantly
finite.\end{thm}

\noindent {\bf Proof.}\quad As we noted above, the proof  here is a
triangulated version of the proof of Gentle-Todorov's theorem.
Assume that $M\in \mathcal{C}$ is an arbitrary object. Take its left
$\mathcal{Y}$-approximation $y_M: M \longrightarrow Y_M$ with
$Y_M\in \mathcal{Y}$. Form a triangle $K \stackrel{k}\longrightarrow
M \stackrel{y_M}\longrightarrow Y_M \longrightarrow K[1]$. Take a
left $\mathcal{X}$-approximation $x_K: K \longrightarrow X_K$ of
$K$.

\vskip 5pt

 Recall from \cite[Appendix]{K} that the octahedral axiom (TR4)
is equivalent to the axioms (TR4') and (TR4''). Hence we have a
commutative diagram of triangles

\[\xymatrix@C=18pt @R=7.5pt{
K \ar[rr]^{k} \ar[dd]^{x_K} & & M \ar[rr]^{y_M} \ar@{.>}[dd]^{z_M}
&& Y_M \ar[rr]
\ar@{=}[dd] && K[1] \ar[dd]^{x_K[1]}\\
& (**) & & & & &\\
X_K \ar@{.>}[rr]^{i_M}       & & Z_M \ar[rr] & & Y_M \ar[rr] &&
X_K[1] }\]

where the square $(**)$ is a \emph{homotopy cartesian square}, that
is, there is a triangle

\begin{align} \label{equation1} K \stackrel{k\choose
x_K}\longrightarrow M\oplus X_K \stackrel{(z_M, \;
-i_M)}\longrightarrow Z_M \dashrightarrow  K[1].\end{align}

Note that $Z_M\in \mathcal{X}\ast \mathcal{Y}$. We claim that the
morphism $z_M: M \longrightarrow Z_M$ is a left $\mathcal{X}\ast
\mathcal{Y}$-approximation of $M$. Then we are done.

\vskip 5pt

To see this, assume that we are given a morphism $f: M
\longrightarrow Z$, and $Z\in \mathcal{X}\ast \mathcal{Y}$, that is,
there is a triangle $X \stackrel{i}\longrightarrow Z
\stackrel{\pi}\longrightarrow Y \longrightarrow X[1]$ with $X\in
\mathcal{X}$ and $Y\in \mathcal{Y}$. Since $y_M$ is a left
$\mathcal{Y}$-approximation, then the composite morphism $\pi\circ
f$ factors through $y_M$, say there is a morphism $c:
Y_M\longrightarrow Y $ such that $\pi\circ f=c\circ y_M$. Hence by
the axiom (TR3), we have a commutative diagram

\[\xymatrix{
K  \ar[r]^{k} \ar@{.>}[d]^{a} & M \ar[r]^{y_M} \ar[d]^{f}    & Y_M
\ar[d]^{c}
\ar[r] & K[1] \ar@{.>}[d]^{a[1]}  \\
X \ar[r]^{i}   &     Z  \ar[r]^{\pi}   & Y \ar[r]   & X[1] }\]

Since $x_K$ is a left $\mathcal{X}$-approximation, the morphism $a$
factors through $x_K$, say we have $a=a'\circ x_K$ with $a': X_K
\longrightarrow X$. Hence $(i\circ a')\circ x_K=f\circ k$, and thus
$(f, \; -i\circ a')\circ {k\choose x_K}=0$. Applying the
cohomological functor ${\rm Hom}_\mathcal{C}(-, Z)$ to the triangle
(\ref{equation1}), we deduce that there is a morphism $h: Z_M
\longrightarrow Z$ such that
$$h\circ (z_M, \; -i_M)=(f, \; -i\circ
a').$$

 In particular, the morphism $f$ factors through $z_M$, as
required. \hfill $\blacksquare$

\vskip 15pt

\section{Applications of Gentle-Todorov's Theorem}

In this section, we apply Gentle-Todorov's theorem  to the
representation theory of artin algebras. We obtain short proofs of a
classical result by Ringel and a recent result by Krause and
Solberg.

\vskip 5pt

Let $A$ be an artin algebra, $A\mbox{-mod}$ the category of finitely
generated left $A$-modules. Dual to the notions of left
approximations and covariantly finite subcategories, we have the
notions of\emph{ right approximations} and \emph{contravariantly
finite subcategories}. A subcategory is called \emph{functorially
finite}, it is both covariantly finite and contravariantly finite.
All these properties are called \emph{homologically finiteness
properties}. \vskip 5pt

 We need more notation. Let $\mathcal{X}\subseteq A\mbox{-mod}$ be a
 subcategory.  Set ${\rm add}\; \mathcal{X}$ to be
its \emph{additive closure}, that is,  the subcategory consisting of
direct summands of modules in $\mathcal{X}$. Note that the
subcategory $\mathcal{X}$ has these homological finiteness
properties if and only if ${\rm add}\; \mathcal{X}$ does. Let $r\geq
1$ and $\mathcal{X}$ a subcategory of $A\mbox{-mod}$. Set
$\mathcal{F}_r (\mathcal{X})=  \mathcal{X}\ast \mathcal{X} \ast
\cdots \ast \mathcal{X}$ (with $r$-copies of $\mathcal{X}$). Hence a
module $M$ lies in $\mathcal{F}_r(\mathcal{X})$ if and only if   $M$
has a filtration of submodules $0=M_0\subseteq M_1 \subseteq M_2
\subseteq \cdots \subseteq M_r=M$ with each factors $M_i/M_{i-1}$ in
$\mathcal{X}$.

\vskip 5pt

Note that the abelian category $A\mbox{-mod}$  has enough projective
and enough injective objects, and thus Gentle-Todorov's theorem and
its dual (on contravariantly finite subcategories) hold. Thus the
following result is immediate.

\vskip 10pt

\begin{cor}\label{corfirst} {\rm (\cite[Corollary 2.8]{SS})} Let $r\geq 1$ and $\mathcal{X}$
 a subcategory of
$A\mbox{\rm -mod}$. Assume that $\mathcal{X}$ is covariantly finite
(resp. contravariantly finite, functorially finite). Then the
subcategories $\mathcal{F}_r(\mathcal{X})$ and ${\rm add}\;
\mathcal{F}_r(\mathcal{X})$ are covariantly finite (resp.
contravariantly finite, functorially finite).
\end{cor}

\vskip 10pt

Recall that a subcategory $\mathcal{X}$ in $A\mbox{-mod}$ is said to
be a \emph{finite subcategory} provided that there is a finite set
of modules $X_1, X_2, \cdots , X_r$ in $\mathcal{X}$ such that each
module in $\mathcal{X}$ is a direct summand of direct sums of copies
of $X_i$'s. Finite subcategories are functorially finite
(\cite[Proposition 4.2]{AS1}). Let $r, n\geq 1$ and
$\mathcal{S}=\{X_1, X_2, \cdots,X_n\}$  a finite set of modules,
denote by $ \mathcal{S}^{\oplus}$ the subcategory consisting of
direct sums of copies of modules in $\mathcal{S}$; for each $n\geq
1$, set
$\mathcal{F}_r(\mathcal{S})=\mathcal{F}_r(\mathcal{S}^\oplus)$. Note
that $ \mathcal{S}^\oplus$ is a finite subcategory and thus a
functorially finite subcategory. The following is a direct
consequence of Corollary \ref{corfirst}.

\vskip 10 pt

\begin{cor}  \label{corsecond} Let  $r, n\geq 1$, and
 $\mathcal{S}=\{X_1, X_2, \cdots X_n\}$ a finite set of  $A$-modules. Then the
 subcategories $\mathcal{F}_r(\mathcal{S})$ and ${\rm add}\; \mathcal{F}_r(\mathcal{S})$ are functorially
 finite.
 \end{cor}

 \vskip 10pt

Let $n\geq 1$ and $\mathcal{S}$ be as above. Ringel introduces in
\cite{R} the subcategory $\mathcal{F}(\mathcal{S})$ to be the
subcategory consisting of modules $M$ with a filtration of
submodules $0=M_0\subseteq M_1 \subseteq M_2 \subseteq \cdots
\subseteq M_r=M$ with $r\geq 1$ and each factor $M_i/M_{i-1}$
belonging to $\mathcal{S}$ . One observes that
$\mathcal{F}(\mathcal{S})=\bigcup_{r\geq 1}
\mathcal{F}_r(\mathcal{S})$. Then we obtain the following classical
result of Ringel with a short proof.

\vskip 10pt

\begin{cor} {\rm ({\rm Ringel}, \cite[Theorem 1]{R} and \cite{R2})} Let $n\geq 1$ and $\mathcal{S}=\{X_1, X_2,
 \cdots, X_n\}$ a finite set of $A$-modules. Assume that ${\rm Ext}_A^1(X_i,
 X_j)=0$ for $i\leq j$. Then the subcategory
 $\mathcal{F}(\mathcal{S})$ is functorially finite.
\end{cor}

\noindent {\bf Proof.}\quad First note the following
\emph{factors-exchanging operation}: let  $M$ be a module with  a
filtration of submodules
$$0=M_0\subseteq M_1 \subseteq  \cdots
\subseteq M_{i-1}\subseteq M_{i} \subseteq M_{i+1} \subseteq \cdots
\subseteq M_r=M,$$ and we assume that ${\rm Ext}_A^1(M_{i+1}/M_{i},
M_i/M_{i-1})=0$ for some $i$, then $M$ has a new filtration of
submodules
$$0=M_0\subseteq M_1 \subseteq \cdots \subseteq M_{i-1}
\subseteq M'_{i} \subseteq M_{i+1} \subseteq \cdots \subseteq
M_r=M$$
exchanging the factors at $i$, that is, $M'_i/M_{i-1}\simeq
M_{i+1}/M_i$ and $M_{i+1}/M'_i\simeq M_{i}/M_{i-1}$. \vskip 5pt

 We claim that $\mathcal{F}(\mathcal{S})=\mathcal{F}_n(\mathcal{S})$.
 Then by Corollary \ref{corsecond} we are done. Let $M\in
 \mathcal{F}(\mathcal{S})$. By iterating the factors-exachanging
 operations, we may assume that the module $M$ has a filtration
  $0=M_0\subseteq M_1 \subseteq M_2 \subseteq \cdots \subseteq M_r=M$
such that there is  a sequence of numbers $1\leq r_1\leq r_2\leq
\cdots \leq r_n=r$  satisfying that the  factors $M_j/M_{j-1}\simeq
X_i$ for ${r_{i-1}+1}\leq j\leq r_i$ (where $r_0=0$). Because of
${\rm Ext}_A^1(X_i, X_i)=0$, we deduce that  $M_{r_i}/M_{r_{i-1}+1}$
is a direct sum of copies of $X_i$ for each $1\leq i\leq n$.
Therefore $M\in \mathcal{F}_n(\mathcal{S})$, as required. \hfill
$\blacksquare$

\vskip 10pt

We will give a short proof to a surprising result recently obtained
by Krause and Solberg \cite{KS}. Note that their proof uses
cotorsion pairs on the category of infinite length modules
essentially, while our proof uses only finite length modules. Recall
that a subcategory $\mathcal{X}$ of $A\mbox{-mod}$ is
\emph{resolving} if it contains all projective modules and it is
closed under extensions, kernels of epimorphims and direct summands
(\cite[p.99]{AB}).

\vskip 10pt

\begin{cor}{\rm (Krause-Solberg, \cite[Corollary 0.3]{KS})}
A  resolving contravarianly finite subcategory of {\rm
$A\mbox{-mod}$} is covariantly finite, and thus functorially finite,
\end{cor}

\noindent{\bf Proof.}\quad  Let $\mathcal{X}\subseteq A\mbox{-mod}$
be a resolving contravariantly finite subcategory. Assume that
$\{S_1, S_2, \cdots, S_n\}$ is the complete set of pairwise
nonisomorphic simple $A$-modules, and take, for each $i$, the
minimal right $\mathcal{X}$-approximation $X_i\longrightarrow S_i$.
Set $\mathcal{S}=\{X_1, X_2, \cdots, X_n\}$. Denote by $J$ the
Jacobson idea of $A$, and assume that $J^r=0$ for some $r\geq 1$.
Thus every modules $M$ has a filtration $0=M_0\subseteq M_1\subseteq
\cdots \subseteq M_r=M$ with semisimple factors. Hence by
\cite[Propostion 3.8]{AR2}, or more precisely by the proof
\cite[Proposition 3.7 and 3.8]{AR2}, we have that $\mathcal{X}={\rm
add}\; \mathcal{F}_n(\mathcal{S})$. By Corollary \ref{corsecond} the
subcategory $\mathcal{X}$ is functorially finite.\hfill
$\blacksquare$

\vskip 10pt

Let us end with an example of functorially finite subcategories.

\begin{exm} {\rm Let $A$ be an artin algebra and  $I$  a two-sided ideal of $A$
such that the quotient algebra $A/I$ is \emph{of representation
finite type}, that is, there are only finitely many isoclasses of
(finitely generated) indecomposable $A/I$-modules. For example, the
Jacobson ideal $J$ satisfies this condition. Let $r\geq 1$, and let
$\mathcal{X}_r$ be the subcategory of $A\mbox{-mod}$ consisting of
modules annihilated by $I^r$. We claim that the subcategory
$\mathcal{X}_r$ is functorially finite in $A\mbox{-mod}$.

\vskip 5pt

To see this, first note that the subcategory $\mathcal{X}_1$ could
be identified with $A/I\mbox{-mod}$, and hence  by assumption
$\mathcal{X}_1$ is a finite subcategory, and thus  functorially
finite in $A\mbox{-mod}$. Then it is a pleasant exercise to check
that $\mathcal{X}_r=\mathcal{F}_r(\mathcal{X}_1)$. By Corollary
\ref{corfirst} we deduce  that the subcategory $\mathcal{X}_r$ is
functorially finite.}

\end{exm}

\vskip 10pt

\noindent{\bf Acknowledgement.}\quad The author would like to thank
Prof. Henning Krause very much, who half a year ago asked him to
give a direct proof to their surprising result \cite[Corollary
0.3]{KS}. Thanks also go to Dr. Yu Ye for pointing out the
references \cite{SS, SS2} and special thanks go to Prof. Apostolos
Beligiannis who kindly pointed out to the author that Theorem 1.1
and the proof were originally due to Gentle and Todorov.

\vskip 5pt

\bibliography{}

\end{document}